\documentclass[a4paper,10pt]{article}

\usepackage{amsfonts}
\usepackage{amsmath}
\usepackage{amssymb}
\usepackage{latexsym}
\def\Rk{\mathbb{ R}^k}
\def\Ck{(0,1)^k}
\def\Rkp{\mathbb{ R}_{+}^k}

\def\be{\begin{equation}}
\def\ee{\end{equation}}

\makeatletter
 
\@addtoreset{equation}{section}
 
\makeatother

\newtheorem{theorem}{Theorem}[section]

\newtheorem{lemma}[theorem]{Lemma}
\newtheorem{remark}[theorem]{Remark}
\newenvironment{Remark}{\begin{remark}\rm}{\end{remark}}

\newcommand{\proof}{\noindent{\sl Proof:}\quad}
\newcommand{\qed}{\hfil$\Box$}
\def\tende#1{\vtop{\ialign{##\crcr\rightarrowfill\crcr
              \noalign{\kern-1pt\nointerlineskip}
              \hskip3.pt${\scriptstyle #1}$\hskip3.pt\crcr}}}

\title{Density-Profile Processes Describing Biological Signaling Networks:  
Almost Sure Convergence to Deterministic Trajectories\\ }

\author{Roberto Fern\'andez$^1$\thanks{\tt
    roberto.fernandez@univ-rouen.fr}, Luiz Renato Fontes$^2$\thanks{\tt
    lrenato@ime.usp.br}, E. Jord\~{a}o Neves$^2$ \thanks{\tt
    neves@ime.usp.br}}

\date{July, 2007}

\begin{document}

\maketitle

\thanks{\small{\begin{center}[1] Laboratoire de
  Math\'ematiques Rapha{\"e}l Salem, UMR 6085 CNRS-Universit\'e de
  Rouen, Avenue de l'Universit\'e, B.P.12, F76801 St Etienne du
  Rouvray, France; \\
{[2] University of S\~{a}o Paulo;
Rua do Mat\~{a}o, 1010 - Cidade Universit\'{a}ria
CEP~05508-090 - S\~{a}o Paulo - SP - Brasil}\end{center}}}


\begin{abstract}

We introduce jump processes in $\Rk$, called {\it density-profile
  process}, to model biological
signaling networks. They describe
the macroscopic evolution of finite-size 
spin-flip models with $k$ types of spins interacting through a
non-reversible Glauber dynamics. We focus on the
the $k$-dimensional empirical-magnetization vector 
in the thermodynamic limit, and prove that, within arbitrary
finite time-intervals, its path converges almost surely to a
deterministic trajectory determined by a first-order (non-linear) differential equation. 
As parameters of the spin-flip dynamics change,
the  associated dynamical system may go through  bifurcations, 
associated to {\it phase transitions} in the statistical
mechanical setting. We present a simple example of spin-flip
stochastic model leading to a dynamical system with {\it Hopf} and {\it pitchfork}
bifurcations; depending on the parameter values, the  
magnetization random path can either
converge to a unique stable fixed point, converge to one of a pair of
stable fixed points, or asymptotically evolve close to a deterministic  orbit in $\Rk$.


\end{abstract}

\section{Motivation and introduction}

The models in this paper are motivated by the mathematical problem
of modeling and 
analyzing biological signaling networks in cancer research.
Understanding how cells manage to respond properly to noisy 
signals from its environment is an important  
research challenge in molecular cell biology, which requires a 
{\it systems  biology}  \cite{kitano} point of view, with emphasis on 
the analysis of 
global aspects of the system, like modularity, robustness and control
mechanisms.

Extra-cellular information is often  transmitted through cell-membrane
receptors activated by {\it ligands}, such as hormones, 
neurotransmitters or  growth factors, which may trigger
complex time-dependent 
cascades of internal cellular biochemical transformations
and lead to quite different cellular responses, like cell-cycle, cell
arrest or {\it cellular suicide} (apoptosis) \cite{cell}.
Signaling pathways  had to survive 
strong selective pressures and therefore 
must contain sophisticated control
mechanisms \cite{ty} in order to avoid inappropriate 
responses which are associated to several diseases, like cancer. 

Several signaling pathways databases are available nowadays 
\cite{stke} which collect knowledge about 
components and their putative interactions. 
Information is usually presented as an oriented graph whose
nodes are the pathway components or group of components 
and whose (oriented) edges indicate some sort of interaction like 
activation or a repression. There is usually no detailed 
information about the biochemical mechanisms associated to these
interactions. In fact a single oriented edge  may actually turn out to  
involve several different processes like  
regulations of gene transcription/translation, protein
transformations (like phosphorilation) or active (non-diffusive) transport 
from one cellular compartment to another.
The mathematical challenge is to propose useful 
models of this complex situation, which lend themselves 
to rigorous analysis but at the same time  
provide experimentally verifiable predictions. 
 
Given the lack of precise information about the biochemical processes
involved in most pathways, specially those of interest in cancer
research,  the usual modeling approach based in 
chemical kinetics \cite{murray}
seems to be unrealistic.   In particular, it involves too many assumptions 
whose validity  can not be checked at present. 
We seek, therefore, an alternative approach to derive the differential
equations that grasp the qualitative continuous-time behavior of those networks. 
In turn, these differential equations provide a general framework 
---based on dynamical systems ideas~\cite{gu}---
to analyze important cellular
behaviors like cell-cycle restriction points \cite{ty2} and pathway
control mechanisms \cite{lauf}.

In this paper we present microscopic models of the following situation.
Suppose $k$ types of
molecules interact in a cellular biochemical
 network and that each one may be in one of
two states, say, {\it active} or {\it inactive}. Suppose further that
there are a large number $N$ of molecules of each type and that
the collective state of the $kN$ molecules evolve like a
continuous-time Markov process. This process is not the
usual  finite-volume stochastic 
spin model~\cite{ligg}  because, in view of our motivation, we 
must allow 
asymmetric interactions among different types of
molecules. 
Indeed, it can happen that the 
activation of a molecule of type A is triggered by the enzymatic effect of
an active
B molecule, while the activation of a molecule of type
B does not depend on those of type A.
In view of our biochemical 
motivation it is also natural to assume that
any molecule may interact which any other.  
These assumptions define spin models what we call 
{\it type-dependent interaction models},  which are fomally
presented in Section~\ref{stoSys}. 

The dynamics of the $k$-dimensional
empirical density vector of these spin models are described by a
the so-called density-profile
process introduced below.  
These are random walk jump-processes  
in $\Rk$ with jumps of size $1/N$, whose expected drift velocity 
$V(x)$ does not depend on $N$.  
The main result of this paper (Theorem \ref{theo})
is the proof that the paths of such a process converge almost surely
to the trajectories of the dynamical system having $V$ as the velocity field.
While similar models where considered in the literature~ \cite{durrett}, the
issue of the almost sure convergence has not been addressed before.  As
we show in an example, the resulting dynamical systems can exhibit a very rich
behavior, including bifurcations.

The main mathematical steps in our convergence proof are the following:
\begin{itemize}
\item[(i)] A graphical construction (Section \ref{sr.gc}) 
that allows a coupled simultaneous construction of density-profile
processes for different $N$.  
\item[(ii)] An auxiliary process $\{\widehat{m}_t^{x^0,N}\}_{t\geq 0}$ (Section \ref{sr.ap})
defined through a simple spin-flip model (independent flips with
time-dependent rates) which shadows the deterministic dynamical
system (Lemma \ref{lr.1}).
\item[(iii)] A coupling between the auxiliary and the density-profile processes
that keep both processes as close to each other as possible (Section \ref{sr.cd}).  
Instants where
they move further apart define a process of \emph{discrepancies}.  Bounds on
the rate of these discrepancies yield our convergence theorem (Theorem \ref{theo}).
\end{itemize} 

\section{Convergence of density profile process to dynamical systems}

In this section we define our basic processes and 
present the main mathematical result of the paper.
In next section we shall realize them through the
macroscopic dynamics of stochastic Ising models.

\subsection{Density-profile processes $\{m^{x^0 N}_t\}_{t\ge 0}$}
\label{dpp}

A density-profile process is a  continuous time 
jump process in the hypercube 
$ {\cal D}_N = (-\frac{1}{N},1+\frac{1}{N})^k.$,
for $k, N \geq 1$. At each jump, a point $x \in  {\cal D}_N$ changes
one of its coordinates by $\pm 1/N$,
with rates that depend smoothly on $x$.  These rates are defined
in the following way. 
We start with two bounded Lipschitz functions
$\lambda, \mu: \Rk\longrightarrow\Rkp$ 
which, in turns, define functions  $f,g: \Rk\longrightarrow\Rkp$ through
the relations
\be
\label{f}
f_i(x) =\left\lbrace
\begin{array}{lcl}
(1-x_i) \lambda_i(x)  & \mbox{  if  }  & 0 \leq x_i \leq 1  \\ 
\lambda_i(0^+) & \mbox{  if  } & x_i \leq 0 \\ 
0 & \mbox{  if  } & x_i \geq 1 
\end{array}
\right.
\ee
and 
\be
\label{g}
g_i(x) =\left\lbrace
\begin{array}{lcl}
x_i \mu_i(x)  & \mbox{  if  }  & 0 \leq x_i \leq 1  \\ 
0 & \mbox{  if  } & x_i \leq 0 \\ 
\mu_i(1^-) & \mbox{  if  } & x_i \geq 1 
\end{array}
\right.
\ee
for $i=1,\ldots,k$.

A density-profile process $\{m^{x^0 N}_t\}_{t \geq 0}$ is 
a random-walk  process in ${\cal D}_N$ which starts at $x^0$ and 
evolves in continuous time through jumps of size $1/N$ along the
coordinate directions.  From each position $x$, the rates for jumps forward or
backwards  along the
coordinate direction $i$ are, respectively, $N f_i(x)$ and  $N g_i(x)$. 
That is, for  $1\leq i \leq k$ and $x \in {\cal D}_N$,
\be
\label{rateup}
 N f_i(x) =  \frac{d}{dt} P\Bigl(m^{x N}_t=x+\frac{e_i}{N} \Bigr)\Bigr|_{t=0}
\ee
and
\be
\label{ratedown}
N g_i(x) =  \frac{d}{dt} P\Bigl(m^{x N}_t=x-\frac{e_i}{N} \Bigr)\Bigr|_{t=0}\;,
\ee
where $e_i$ denotes the unit vector
along direction $i$.  

In our applications, the variable $x_1,\ldots,x_k$ represent the densities
of $k$ types of objects that can be present at $N$ different sites.  The function
$\lambda_i$ is the rate of creation or activation  
of an object of type $i$ at a site where the object is absent or inhibited.
The function $\mu_i$ is the rate for the opposite move.  Thus, $Nf_i$ 
(resp.\ $g_i$) is the total rate of creation (resp.\ destruction) of an
object of type $i$ throughout the entire collection of sites.  

In the sequel we assume the choice of a common probability space where the
density-profile processes (and other auxiliary processes defined below) 
are simultaneously realized for all $N$ and $x^0$.
The corresponding probability measure will be denoted $P$.  
(The graphical
constructions introduced in Section \ref{sr.pct} offer, in fact, a concrete way of defining
this common space.)

\subsection{Convergence to a dynamical system $\{x^{x^0}_t\}_{t\ge 0}$}
\label{dynSys}

Let $\{m^{x^0 N}_t\}_{t \geq 0}$ be the density-profile  process in
${\cal D}_N$ defined as above for a given pair of appropriate functions
$(\lambda,\mu)$, and let $V: \Rk\longrightarrow \Rkp$ be its associated
{\it velocity field}:
\be
\label{vel}
V(x)=\lim_{t\downarrow 0} \frac{{\bf E}(m^{x N}_t-x)}{t}=f(x)-g(x)\;,
\ee
that is,
\be
\label{V}
V_i(x) =
\left\lbrace
\begin{array}{lcl}
(1-x_i) \lambda_i(x) - x_i \mu_i(x) & \mbox{  if  }  & 0 \leq x_i \leq 1  \\ 
\lambda_i(0^+) & \mbox{  if  } & x_i \leq 0 \\ 
\mu_i(1^-) & \mbox{  if  } & x_i \geq 1 
\end{array}
\right.
\ee
Let $\{x^{x^0}_t\}_{t \geq 0}$ be the solution of the dynamical system
\be
\label{ds}
\dot{x}_t=V(x_t)
\ee
starting at $x^0 \in \Ck$. 
The global trajectory exists by the smoothness of
the field $V$.  Furthermore, the flow does not leave $\Ck$ because
$ V_i(0^+) >0$ and $V_i(1^-)<0$
for $i =1, \ldots k$.

The main result of this paper is the convergence of the sequence of
density profile processes $(m^{x^0 N}_t)_N$ to the trajectory $x^{x^0}_t$.
For $\epsilon > 0$ let ${\tau}^N_{\epsilon}$ be the stopping time 
\be
\label{tau}
{\tau}^N_{\epsilon}= \inf\Bigl\{ t \geq 0:|m_t^{x^0,N}-x_t^{x^0}| > 
\frac{1}{N^{  \frac{1}{2}-\epsilon }}\Bigr\} 
\ee
and, for a given $T$, $0\leq T<\infty$, 
write ${\cal A}^{T}_{N \epsilon}=\{\tau^N_{\epsilon} < T\}$.
The following is our main result.

\begin{theorem} \label{theo} 
Let $\lambda$ and $\mu$ be bounded functions from $\Rk$ to
$\Rkp$ satisfying the Lipschitz condition 
\be
\label{lipchitz}
|\lambda(x)-\lambda(y)|\leq K |x-y| \quad;\quad
|\mu(x)-\mu(y)|\leq K |x-y| 
\ee
for some $K>0$ and all $x,y\in (0,1)^k$.
Then, for any finite $T$, initial position $x^0$ and $\epsilon >0$,
\be
\label{er:theo}
P\Bigl({\overline{\lim}}_{N} {\cal A}^{T}_{N \epsilon}\Bigr)=0
\ee
\end{theorem}

That is, for typical realizations there exists some $N_{\epsilon,T}$ such that 
for $N>N_{\epsilon,T}$ each process $\{m^{x^0 N}_t\}_{t \geq 0}$ stays
within a distance $N^{-1/2+\epsilon}$ of the deterministic path $\{x_t^{x^0}\}_{t \geq 0}$
at least up to time $T$.  

Dynamical systems of the form (\ref{ds})/(\ref{V}) 
can exhibit quite complex dynamics 
---even for simple choices of $\lambda$ and $\mu$---, 
including stable orbits and chaotic
behavior. In Section~\ref{cim} we present an example of a
system leading to a Hopf bifurcation.

\section{Associated stochastic spin models}
\label{stoSys}

We present here a family of stochastic spin models, motivated by biological cellular systems,
whose empirical densities evolve as density-profile processes.

\subsection{General definition}
We consider a discrete set $\Lambda$ of \emph{sites} and a finite set 
${\cal T}=\{1,\ldots, k\}$ of \emph{types}.   Each type $i\in\mathcal{T}$ can be present 
at each {\it site} $\ell\in \Lambda$.  Thus, we choose our configuration space as
$\Sigma=\{-1,+1\}^{{\cal T}\times \Lambda}$, where for each $\eta \in
\Sigma$, the value $\eta(i,\ell)=+1$ ($\eta(i,\ell)=-1$) indicates the presence (absence) 
of a particle of type $i$ at site $\ell$. 

We consider continuous-time processes in $\Sigma$
in which only single spin flips are allowed. 
For $\eta \in \Lambda$ and $(i,\ell)\in {\cal T} \times \Lambda$, 
let $\eta^{(i,\ell)}$ denote the  
configuration which is equal to 
$\eta$ except at site $\ell$ where the spin of type $i$ is {\it flipped}.
\be
\label{flip}
\eta^{(i,\ell)}(j,n) =\left\lbrace
\begin{array}{rl}
\eta(j,n)  & \mbox{  if  }  (j,n)\not= (i,\ell)   \\ 
-\eta(i,\ell)  & \mbox{  otherwise  }   
\end{array}
\right.
\ee
We shall consider {\it type-dependent interaction models}, for which
the flip rate for the spin-flip transition 
$\eta \rightarrow \eta^{(i,\ell)}$ is a monotone
non-increasing function of 

\be
\label{delta}
\Delta(i,\ell,\eta)=H_i(\eta^{(i,\ell)})-H_i(\eta)
\ee

\noindent for a  {\it Hamiltonian vector} $H(\eta)= (H_1(\eta),
\ldots, H_k(\eta))$.  

A simple choice is 

\be
\label{rate}
c(i,\ell,\eta) = \exp\Bigl[- \Delta(i,\ell,\eta)  \Bigr].
\ee

We focus on the particular case
\be
\label{h}
H_i(\eta)=  - \frac{1}{2} \sum_{\ell\in \Lambda} \Biggl (\sum_{(j,n)\in {\cal
    T}\times \Lambda} \frac{ \alpha_{j i}}{|\Lambda|}\, \eta(j,n)\, \eta(i,\ell) 
+ a_i \,\eta(i,\ell) \Biggr )\;.
\ee
This corresponds to a {\it mean-field} interaction, where
$\alpha_{j i}$ is the strength of the influence of spins of type $j$
on those of type $i$ and each $a_i$ acts as a type-dependent external
field.  The most interesting phenomena appear when
$\alpha$ is not symmetric. 
The {\it empirical density profile} of a configuration $\eta$, 
is the $k$-uple $m(\eta)=(m_1(\eta), \ldots , m_k(\eta)) \in
\Rkp$, where 
\be
\label{rho}
m_i(\eta)= \frac{|\{\ell\in\Lambda: \eta(i,\ell)=+1\}|}{|\Lambda|}.
\ee
for $1\leq i \leq k$.

Let $\{\sigma^{\eta_0}_t\}_{t \geq 0}$ be the spin-flip process starting
at the configuration $\eta_0 \in \Sigma$ with rates defined by
(\ref{rates}) and (\ref{h}). 
  Then, if $|\Lambda|=N$,
the density-profile process $m(\sigma^{\eta^0}_t)$  
approximates, in the sense of Theorem \ref{theo} the dynamical system
(\ref{ds})/(\ref{V}), defined by
\be 
\label{llambda}
\lambda_i(x)=\exp\biggl(\sum_{j\in {\cal T}} \alpha_{j i} x_j +a_i\biggr)
\ee
and
\be 
\label{mmu}
\mu_i(x)=\exp\biggl(-\sum_{j\in {\cal T}} \alpha_{j i} x_j -a_i\biggr).
\ee

Given our biological motivation, models with mean-field interaction
like (\ref{h}) are natural. They are also the simplest models with
type-dependent rates of the form (\ref{rate}). From the mathematical point of
view, though, the analysis of 
models with local interaction \cite{ruelle}, possibly with {\it
  stirring} \cite{pablo} are much more interesting.

\subsection{Example: Cyclic-interaction model}
\label{cim}

We exhibit now a model ---called the {\it cyclic-interaction model}---
defined through a
simple choice of the interaction matrix $\alpha$ which nevertheless leads to
interesting (deterministic) dynamical behavior.
Think $\{1,\ldots, k\}$ as points on the circle and, for each $i\in\mathcal{T}$ let $c(i)$ denote
the nearest-neighbor of $i$ in the counter-clockwise direction.
We assume that  $\alpha_{j i}= 0$ unless $j=c(i)$ and that, furthermore, 
all non-zero terms in $\alpha$ have the same absolute value.  That is,
\be 
\label{cyclic}
\alpha_{j i}  =\left\lbrace
\begin{array}{cl}
s_i\,J   & \mbox{  if  } j=c(i)  \\ 
0   & \mbox{  otherwise   } 
\end{array}
\right.
\ee
where $s_i \in \{-1,+1\}$ representing the signals, and $J>0$.   We also set 
$a_i=-J/2$, for $1\leq i \leq k$.   In this way,  once the signs
$\{s_i\}_{i=1}^{k}$ are chosen,  $J$ is the only free parameter of the model.  

If $s_i=1$, the rate with which spins of type $i$ flip from $-1$ to $+1$ [defined in
(\ref{llambda})], is an increasing  function of $x_{c(i)}$, 
the density of spins $+1$ of type $c(i)$. Borrowing 
statistical mechanical nomenclature, we say that the interaction of spins of
type $c(i)$ with those of type $i$ is {\it ferromagnetic} \cite{thom}. In the
biochemical context, where $x_i$ measures the density of
some chemical component $i$, this means
that the component $c(i)$ {\it activates} the production of component $i$.
On the other hand, if $s_i=-1$ the rate
for a spin of type $i$ to flip from $+1$ to $-1$ [defined in
(\ref{mmu})], decreases as a function of $x_{c(i)}$, and 
 the interaction of spins of
type $c(i)$ with those of type $i$ is {\it anti-ferromagnetic}. In
biochemical terms, the component $c(i)$ inhibits the
production of component $i$. 

The dynamical system (\ref{ds}) associated
to the cyclic-interaction model (\ref{cyclic}) is:
\be
\label{dsc}
\dot{x}_i= e^{s_i J(x_{c(i)}-\frac{1}{2})} -x_i 
\Bigl(e^{s_i J(x_{c(i)}-\frac{1}{2})}+e^{-s_i J(x_{c(i)}-\frac{1}{2})}\Bigr)
\ee
for $1\leq i \leq k$. 
If $J$ is small, this system has a single stable equilibrium  point at 
$(1/2, \ldots , 1/2)\in \Rk$, whichever the choice of signs $s_i$.
For larger $J$, the behavior of the dynamical system (\ref{dsc})
crucially depends on whether the product of signals is positive or negative.  
If $\Pi_{i=1}^{k} s_i =-1$ ---a {\it frustrated} model  in statistical mechanical terms---
there is no (global) density-profile where all pairs of types of spins 
minimize their mutual interaction.  This system exhibits a Hopf
bifurcation \cite{gu} as $J$ exceeds a critical value, which depends on $k$. 
In the non-frustrated case, the model behaves as the {\it Curie-Weiss} model. 
Formally:
\begin{theorem}
Consider the dynamical system (\ref{dsc}) with $k\geq 3$
\begin{itemize}
\item[(a)] If  $\Pi_{i=1}^{k} s_i=1$, there is a 
bifurcation at $J_c=2$: the fixed point
$(1/2,\ldots, 1/2)$ looses stability and two stable points appear for
$J> J_c$. 
\item[(b)]  If  $\Pi_{i=1}^{k} s_i=-1$, there is a {\it Hopf} bifurcation at
$J_c=2/\cos(\pi/k)$.
\end{itemize}
\end{theorem}

\proof Write  $s=\Pi_{i=1}^{k} s_i$.
A simple computation shows that 
near  ${\bf 1/2}= (1/2, \ldots, 1/2)\in \Rk$ the dynamical
system (\ref{dsc}) is close to $\dot{x}=A (x-{\bf 1/2})$, 
where $A$ is a $k\times k$ real
matrix with eigenvalues $s J e^{\frac{2\pi l}{k}i}-2$, $l=0,
1, \ldots ,k-1$. The fixed point is stable if the real parts of all eigenvalues are
strictly negative, and stability is lost when one one of the real parts
becomes positive. Thus, if $s=1$ the fixed point ${\bf 1/2}$ loses
stability at $J_c=2$ when the eigenvalue corresponding to $l=0$
crosses the imaginary axis through the origin. On the other hand, if
$s=-1$, the stability is lost when two
eigenvalues, symmetric around the real axis,  cross the imaginary
axis. This occurs at $J_c=2/\cos(\pi/k)$. \qed

\begin{Remark} For instance, if $k=3$ and 
all interactions are antiferromagnetic ($s_i=-1$ for $i=1,2,3$), the
dynamical system has stable orbits for $J>J_c=4$. The convergence result,
Theorem \ref{theo}, implies that, within any finite time interval,
the density-profile process evolves as close to this orbit  as wished, for $N$
sufficiently large.
\end{Remark}

\section{Proof of the convergence theorem}\label{sr.pct}
\subsection{The auxiliary process $\{\widehat{m}_t^{x^0,N}\}_{t\geq 0}$}
\label{sr.ap}

To prove Theorem \ref{theo} we introduce an auxiliary stochastic spin
model with independent spins flips but {\it time-dependent rates}.

Let  $\Lambda=\{1,\ldots,N\}$ and let
$\{\eta_t(i,n) : (i,n) \in \mathcal{T}\times\Lambda\}_{t\geq 0}$, 
 be $kN$ independent Markov chains with
state space $\{-1,+1\}$. Thus, for each $t\geq 0$, $\eta_t \in \Sigma
= \{-1,+1\}^{{\cal T} \times \Lambda}$, with ${\cal
  T}=\{1,\ldots,k\}$, as defined in the previous section.
For each Markov chain $\{\eta_t(i,n)\}_{t\geq
   0}$ the flips from $-1$ to $+1$ and from $+1$ to $-1$ 
 have time-dependent rates given,
 respectively, by $\lambda_i(x^{x^0}_t)$ and $\mu_i(x^{x^0}_t)$,
where $\{x^{x^0}_t\}_{t\geq 0}$ is the solution of the dynamical system
(\ref{ds}) from the initial position $x^0$.  We initialize these chains
with the uniform distribution on 
configurations $\eta_0$ with $m(\eta_0)=x^0$, where 
$m(\eta_t)$ is defined in (\ref{rho}).  
The total number of spins of each type $i$ is, thus, fixed 
and equal to $x^0_i$, but he initial density components
 $m_1(\eta_0), \ldots, m_k(\eta_0)$ are independent.
 We denote $\{\widehat{m}_t^{x^0,N}\}_{t\geq 0}$ the corresponding
 density-profile process.  
 
We observe that the chain at each $(i,n)\in\mathcal{T}\times\Lambda$
satisfies Kolmogorov's equation.  Hence, for
 $p_t(i,n)=P(\eta_t(i,n)=+1)$,  we have 
\be
\label{kol}
\dot{p}_t(i,n)\;=\;[1-p_t(i,n)]\,\lambda_i(x_t^{x^0}) - p_t(i,n)\,
\mu_i(x_t^{x^0})\;.
\ee
Therefore each function 
$t\to p_t(i,n)$ is a solution of the differential equation
(\ref{ds}) with $V$ as in (\ref{V}).  Hence 
\be \label{er:2}
p_0(i,n)=(x^0)_i 
\quad\Longrightarrow\quad 
p_t(i,n)=(x_t^{x^0})_i \quad \forall t\ge 0\;,
\ee
for all $i \in \mathcal{T}$ and $n\in\Lambda$
[$(y)_i$ indicates the $i$-th component of vector $y \in \Rk$].
While (\ref{er:2}) is true for the auxiliary process $\{\widehat{m}_t^{x^0,N}\}$,
we are interested in following the actual empirical densities.  
Next lemma proves that also the path followed by these densities 
remain, at all times, close to the trajectories of the dynamical system.
\begin{lemma} \label{lr.1}For $\delta >0$ there exists $c>0$ such that 
\be
\label{bound2}
P\Bigl(\bigl|\widehat{m}_t^{x^0,N}-x^{x^0}_t\bigr|>
N^{\delta-1/2}\Bigr)\; <\; \exp (-c N^{\delta})
\ee
for $t\geq 0$.
\end{lemma}

\proof
Let us introduce yet another auxiliary process, denoted 
$\{\widehat{m}_t^{b(x^0),N}\}_{t\geq 0}$, defined exactly as
$\{\widehat{m}_t^{x^0,N}\}_{t\geq 0}$ but with initial spins
chosen independently with $P(\eta_0(i,n)=+1)=(x^0)_i$.
Hence, the density components $m_1(\eta), \ldots, m_k(\eta)$ 
are independent and each $m_i(\eta)$ 
has a binomial distribution
with parameters $N$ and $p_i=(x^0)_i$,  $i=1,\ldots, k$. 
(This means that, for large $N$,  
$\widehat{m}_t^{b(x^0),N}$ starts at a random
position in ${\cal D}$ close to $x^0$, while 
$\widehat{m}_t^{x^0,N}$ starts precisely at $x^0$.)

This new auxiliary process also satisfies  (\ref{er:2}) but has the
advantage that the spin chains remain independent at all times,
and, by (\ref{er:2}), the proportions of spins of each type coincide
with the components of $x^{x^0}_t$.  Therefore
$N \widehat{m}_t^{b(x^0),N}$ is a vector of independent binomial random variables 
\be
\label{binomial}
N \Bigl(\widehat{m}_t^{b(x^0),N}\Bigr)_i\sim\; {\rm Bin}(N,(x_t^{x^0})_i)
\ee
for $t \geq 0$ and $i\in\mathcal{T}$.  
In particular the the variances of the components of
$\widehat{m}_t^{b(x^0),N}$ are proportional to $1/N$.  Thus, 
the large-deviation properties of binomial distributions \cite{ellis}
imply that for any $\delta>0$ there exists a constant $c$ such that
\be
\label{bound}
P\Bigl(\bigl|\widehat{m}_t^{b(x^0),N}-x^{x^0}_t\bigr|>
\frac12 N^{\delta-1/2}\Bigr) \;<\; \exp (-c N^{\delta})
\ee
for any $t\geq 0$. 

To conclude the proof of the lemma we must show that
both auxiliary processes $\widehat{m}_t^{b(x^0),N}$ and
$\widehat{m}_t^{x^0,N}$ remain close to each other. 
This is more easily done through a coupling argument~\cite{tor,ligg}.
We construct a coupled realization $(\eta^{b(x^0,N)}_t,\eta^{x^0,N}_t)$
of the spin systems defining both processes as follows.  
Spins in both systems flip  
with the same time-dependent rates given in (\ref{kol}).
At sites $(i,n)$ with $\eta_0^{b(x^0),N}(i,n)=\eta_0^{x^0,N}(i,n)$, the spins
evolve together.  
Otherwise, the spins for both processes evolve independently 
until one of them makes a transition, thereby bringing them to a 
common value.  They evolve together ever after.  
As the distance between the corresponding coupled density profiles decreases with time,
\be
\label{1}
\begin{split}
\Bigl| m\bigl(\eta_t^{b(x^0),N}\bigr)-m\bigl(\eta_t^{x^0,N}\bigr)\Bigr| &\leq
\Bigl |m\bigr(\eta_0^{b(x^0),N}\bigr)-m\bigl(\eta_0^{x^0,N}\bigr)\Bigr|\\
&=\Bigl| m\bigl(\eta_0^{b(x^0),N}\bigr)-x^0\Bigr|
\end{split}
\ee
and, therefore,
\begin{eqnarray}
\label{dif}
\lefteqn{
P\Bigl(\bigl |\widehat{m}_t^{x^0,N}-x^{x^0}_t \bigr|>
N^{\delta-1/2}\Bigr) \;\leq }\nonumber\\
&&P\Bigl(\bigl|\widehat{m}_t^{x^0,N}-\widehat{m}_t^{b(x^0),N}\bigr|>
\frac{1}{2}N^{\delta-1/2}\Bigr)+
P\Bigl(\bigl|\widehat{m}_t^{b(x^0),N}-x^{x^0}_t\bigr|>
\frac{1}{2}N^{\delta-1/2}\Bigr)\;.\nonumber\\
\end{eqnarray}
To prove (\ref{bound}) we bound the right-hand side 
using  (\ref{1}) and (\ref{bound}) (for $t=0$) for the first term
and again (\ref{bound}) for the second one.
\qed
\medskip

To prove Theorem \ref{theo} we will show  that, for $N$ large, 
$\widehat{m}_t^{x^0,N}$ and ${m}_t^{x^0,N}$ remain close within arbitrary
finite time intervals. To achieve this we will couple
both stochastic evolutions with the help of a  graphical
construction.

\subsection{Graphical construction: The process $\{g_t^{x^0,N}\}_{t\geq 0}$} \label{sr.gc}
We resort to a graphical construction of  density-profile processes with different $N$
through time-rescaling of auxiliary processes 
$\{g_t^{x^0,N}\}_{t\geq 0}$.  The latter is defined through paths determined by 
Poissonian ``marks''. 
This construction will be adapted in next
section to couple the processes $\widehat{m}_t^{x^0,N}$ and ${m}_t^{x^0,N}$.

To each $y \in {\cal D}_N$ we associate $2k$ independent Poisson
processes: $N_t^{1+}(y)$, $N_t^{1-}(y)$, $\ldots$, $N_t^{k+}(y)$,
$N_t^{k-}(y)$, where each $N_t^{i+}(y)$ has rate $f_i(y)$ and each $N_t^{i-}(y)$ rate
$g_i(y)$. We associate a particular type of mark for the events of 
each type of process and place these marks along the time axis of $y$. 
A Poisson mark associated to the process $N_t^{i+}(y)$
($N_t^{i-}(y)$) carries the instruction to
jump along the positive (negative) $i$ coordinate direction.  

The process $\{g_t^{x^0,N}\}_{t\geq0}$ is defined by 
{\it open paths} in ${\cal D}_N \times \mathbb{R}_+$ determined by the marks. 
These are piecewise linear curves that  
move along the positive time axis until a Poisson mark is
met.  At these times the trajectory moves by $\pm1/N$ along
a coordinate direction according to the type of mark.  
The process $\{g_t^{x^0,N}\}_{t\geq 0}$ is at
position $x$ at time $t$ if   there
exists an open path from $(x^0,0)$ to $(x,t)$.

We see that the evolution of $\{g_t^{x^0,N}\}$
differs from that of $\{m_t^{x^0,N}\}$ only in
that the rates of the latter [given by (\ref{rateup}) and (\ref{ratedown})]
are $N$ times those of the former. 
Thus, one process can be constructed from the other by 
a simple change in the time scale:
\be
m_t^{x^0,N}= g_{Nt}^{x^0,N}\;.
\ee
In words: a \emph{density-profile time} $t$ corresponds to a 
\emph{graphical-construction time} $Nt$.

\subsection{Main coupling and discrepancy process}

We now use the graphical-construction strategy
to produce coupled realizations of the density-profile processes 
${m}_t^{x^0,N}$ and $\widehat{m}_t^{x^0,N}$ 
with an appropriate control of their distance. 
Our coupling forces both processes to keep their relative distance
as much as possible, evolving as a rigid system. 
Of course, since their rates are not equal, they will
make occasional asynchronous moves that may 
take them increasingly apart with the passing of time.    
The coupling is designed so to minimize this asynchrony. 

The coupling involves a number of Poissonian mark processes at different
sites which are updated every time there is an asynchronous move. 
The successive times of these moves correspond to a sequence
of stopping times $\{\tau_n\}_{n  \geq 1}$; the coupling is defined in a  recursive
fashion within successive time intervals $[\tau_{n-1}, \tau_n)$, $n
\geq 1$.
The auxiliary processes, which arise directly from such graphical
coupled construction will be denoted, respectively, by 
${g}_t^{x^0,N}$ and $\widehat{g}_t^{x^0,N}$.  They differ from the
density profiles ${m}_t^{x^0,N}$ and $\widehat{m}_t^{x^0,N}$
only in the time scale, which in the graphical construction is slower by a factor $N$.

\paragraph{Initial stage of the coupling}

Initially, $g_0^{x^0}=\widehat{g}_0^{x^0}=x^0$ and up to the first
stopping time $\tau_1$ (to be defined) we couple them through 
what is known as {\it basic coupling} in particle systems.
For each $y \in {\cal D}_N$ and coordinate direction
$i=1, \ldots, k$ we define six Poissonian mark processes:
\begin{itemize}
\item[(i)] \emph{Marks associated to jumps from $y$ to $y+e_i/N$}:
They are defined by independent Poisson processes
$\widehat{N}_t^{i,+}(y)$, $\widehat{E}_t^{i,+, m}(y)$ and $\widehat{E}_t^{i,+,
  \widehat{m}}(y)$ with respective rates
\be
\begin{split} 
\widehat{u}_t^{i,+}(y) &= \min \bigl\{(1-y_i)\,
  \lambda_i(y)\,,\, (1-y_i)\, \lambda(x_{t/N}^{x^0})\bigr\}\;,\\
\widehat{e}_t^{i,+, m}(y) &= \bigl|(1-y_i) \,\lambda_i(y) - 
\widehat{u}_t^{i,+}\bigr|_+\qquad \mbox{and}\\    
\widehat{e}_t^{i,+, \widehat{m}}(y) &= \bigl|(1-y_i)\, \lambda_i(x_{t/N}^{x^0}) -
\widehat{u}_t^{i,+}\bigr|_+\;.  
\end{split}
\ee
[$\left|z\right|_+=\max\{z,0\}$].
Note the rescaling in time for the deterministic path $\{x_t^{x^0}\}$
needed to represent it on the graphical construction time scale. 

\item[(ii)] \emph{Marks associated to jumps from $y$ to $y-e_i/N$}:
Defined by three independent Poissonian processes 
$\widehat{N}_t^{i,-}(y)$, $\widehat{E}_t^{i,-, m}(y)$ and $\widehat{E}_t^{i,-,
\widehat{m}}(y) $ which are independent from the precedent ones and
have respective rates
\be
\begin{split}
\widehat{u}_t^{i,-}(y) &= \min\bigl\{y_i\,
  \mu_i(y)\,,\, y_i \,\mu(x_{t/N}^{x^0})\bigr\},\\
\widehat{e}_t^{i,-, m}(y) &= \bigl|y_i \,\mu_i(y) - \widehat{u}_t^{i,-}
\bigr|_+\qquad \mbox{and}\\   
\widehat{e}_t^{i,-, \widehat{m}}(y) &= \bigl|y_i \,\mu_i(x_{t/N}^{x^0}) -
\widehat{u}_t^{i,-}\bigr|_+\;.  
\end{split}
\ee

\end{itemize}
As before, we think that occurrence of each of these processes are
associated to particular marks indicating where to jump. 
The jumps of the process $\{g_t^{x^0,N}\}$ occur at the marks of 
$\{\widehat{E}_t^{i,+, m}(y)\}$
  and $\{\widehat{E}_t^{i,-, m}(y)\}$; while those of the process
$\{\widehat{g}_t^{x^0,N}\}$ are at  
$\{\widehat{E}_t^{i,+, \widehat{m}}(y)\}$ and 
$\{\widehat{E}_t^{i,-,  \widehat{m}}(y)\}$.
The marks of the four processes $\{\widehat{E}_t^{i,+, m}(y)\}$,
 $\{\widehat{E}_t^{i,+, \widehat{m}}(y)\}$,   
$\{\widehat{E}_t^{i,-, m}(y)\}$ and  
 $\{\widehat{E}_t^{i,-, \widehat{m}}(y)\}$
 are thus seen by only one of $\{g_t^{x^0,N}\}$ or $\{\widehat{g}_t^{x^0,N}\}$
 and will be called {\it discrepancies}.  The \emph{basic Poisson processes} 
 $\{\widehat{N}_t^{i,\pm}(y)\}$, on the other hand, 
are introduced to ensure that 
$\{g_t^{x^0,N}\}$ and $\{\widehat{g}_t^{x^0,N}\}$ remain equal until they find
the first discrepancy.  This defines a stopping time $\tau_1$ at which
the processes get separated by a distance of $1/N$.  At this time
we can not continue using  the basic coupling.

Formally, we define a \emph{first-discrepancy process}
\be
\label{er.fd}
D^{0}_t \;=\; \sum_{i=1}^k \Bigl[\widehat{E}_t^{i,+, m}(y^0) + 
\widehat{E}_t^{i,-, m}(y^0) + \widehat{E}_t^{i,+, \widehat{m}}(y^0) +
\widehat{E}_t^{i,-,  \widehat{m}}(y^0)\Bigr]
\ee
where $y^0$ is the density-profile path defined by the preceding (level-0) construction.
The \emph{first discrepancy time} $\tau_1$ is the time of the first event of
this process.  A new coupling definition must be introduced at this time, which will be
applied until the second discrepancy time $\tau_2$.  This iterative procedure 
continues up to the time $T$ chosen in Theorem \ref{theo}.  We now present
the recursive step in the definition of such a coupling.

\paragraph{$l$-th stage of the coupling}

Suppose  that the graphical construction has been defined up to 
time $\tau_l$, $l \geq 1$, determining $x^l,\Delta^l\in\mathcal{D}_N$ such that 
\be
\label{er.lt}
g_{\tau_l}^{x^0,N}=x^l \quad,\quad  
\widehat{g}_{\tau_l}^{x^0,N}= x^l + \Delta^l\;.
\ee
[Thus, 
$m_{\tau_l/N}^{x^0,N}=x^l$ and $\widehat{m}_{\tau_l/N}^{x^0,N}=x^l+\Delta^l$.]
From time $\tau_l$ we start a new 
graphical construction, which defines the evolution of both
processes until the next discrepancy appears at time $\tau_{l+1}$. 
We define the following Poisson mark processes for each 
$y \in {\cal D}_N$ and coordinate direction $i=1, \ldots, k$:
\begin{itemize}
\item[(i)] \emph{Marks associated to jumps from $y$ to $y+e_i/N$}:
Let $\widehat{N}_t^{i,+, m}(y)$, $\widehat{N}_t^{i,+, \widehat{m}}(y)$,
 $\widehat{E}_t^{i,+, m}(y)$ and $\widehat{E}_t^{i,+, \widehat{m}}(y)$ 
be  Poisson processes with respective rates
\be
\label{+}
\begin{split} 
\widehat{u}_t^{i,+, m}(y,\Delta^l) & = \min \bigl\{(1-y_i)\,
  \lambda_i(y)\,,\, (1-y_i-\Delta_i^l) \,\lambda(x_{t/N}^{x^0})\bigr\}\;,\\
\widehat{u}_t^{i,+, \widehat{m}}(y,\Delta^l) &= \min\bigl\{(1-y_i)\, 
  \lambda_i(x_{t/N}^{x^0})\,,\, (1-y_i+\Delta_i^l)\, \lambda_i(y-\Delta^l)\bigr\}\;,\\ 
\widehat{e}_t^{i,+, m}(y,\Delta^l) &= \bigl|(1-y_i) \,\lambda_i(y) -
\widehat{u}_t^{i,+,m}(y,\Delta^l)\bigr|_+\qquad \mbox{and}\\   
\widehat{e}_t^{i,+, \widehat{m}}(y,\Delta^l)&= \bigl|(1-y_i) \,\lambda_i(x_{t/N}^{x^0}) -
\widehat{u}_t^{i,+,\widehat{m}}(y,\Delta^l)\bigr|_+\;.
\end{split}
\ee
We observe that $\widehat{u}_t^{i,+, m}(y,\Delta^l) = \widehat{u}_t^{i,+
  ,\widehat{m}}(y+\Delta^l,\Delta^l)$ for any $y \in {\cal D}$, so we identify
\be
\label{er:i1}
\widehat{N}_t^{i,+, m}(y) = \widehat{N}_t^{i,+ ,\widehat{m}}(y+\Delta^l)\;.
\ee
Except for this identification, the different processes are mutually independent and
independent of all previous Poisson mark processes.

\item[(ii)] \emph{Marks associated to jumps from $y$ to $y-e_i/N$}:
They are determined by Poisson processes 
$\widehat{N}_t^{i,-, m}(y)$, $\widehat{N}_t^{i,- ,\widehat{m}}(y)$, 
 $\widehat{E}_t^{i,-, m}(y)$ and $\widehat{E}_t^{i,-, \widehat{m}}(y)$, 
respectively with rates
\be
\label{-}
\begin{split}
\widehat{u}_t^{i,-, m}(y,\Delta^l) &= \min\bigl\{y_i
  \mu_i(y)\,,\, (y_i+\Delta_i^l) \,\mu(x_{t/N}^{x^0})\bigr\}\;,\\ 
\widehat{u}_t^{i,-, \widehat{m}}(y,\Delta^l) &= \min\bigl\{y_i\,
  \mu_i(x_{t/N}^{x^0})\,,\, (y_i-\Delta_i^l) \,\mu(y-\Delta^l)\bigr\}\;,\\ 
\widehat{e}_t^{i,-, m}(y,\Delta^l) &= \bigl|y_i\, \mu_i(y) -
\widehat{u}_t^{i,-,m}(y,\Delta^l)\bigr|_+\qquad \mbox{and}\\   
\widehat{e}_t^{i,-, \widehat{m}}(y,\Delta^l) &= \bigl|y_i \,
\mu_i(x_{t/N}^{x^0}) - \widehat{u}_t^{i,-,\widehat{m}}(y,\Delta^l)\bigr|_+\;; 
\end{split}
\ee
with the identification
\be
\label{er:i2}
\widehat{N}_t^{i,-, m}(y) = \widehat{N}_t^{i,- ,\widehat{m}}(y+\Delta^l)\;.
\ee
All these processes are independent among themselves, except for the
preceding identification, and independent of other mark processes.
\end{itemize}

The process $\{g_t^{x^0,N}\}$ jumps only at the marks 
placed by the processes $\{\widehat{E}_t^{i,\pm, m}(y)\}$, while
process $\{\widehat{g}_t^{x^0,N}\}$ does so at the marks of
$\{\widehat{E}_t^{i,\pm,\widehat{m}}(y)\}$.  
Due to identifications (\ref{er:i1})/(\ref{er:i2}), the basic Poisson marks
$\{\widehat{N}_t^{i,\pm, m}(y)\}$ seen by 
$\{g_t^{x^0,N}\}$ at a given position $y$ coincide with the basic marks
seen by $\{\widehat{g}_t^{x^0,N}\}$ at its corresponding 
position $y+\Delta^l$.  Hence, the two graphic processes evolve rigidly, 
keeping a separation $\Delta_l$, until a discrepancy is met, that is, 
until one of the processes
responds to a Poisson mark that the other ignores. This happens either
because $\{g_t^{x^0,N}\}$ at a certain position $y$ meets a
mark of $\{\widehat{E}_t^{i,+, m}(y)+ \widehat{E}_t^{i,-, m}(y)\}$ or because
$\{\widehat{g}_t^{x^0,N}\}$, at the corresponding position
$y+\Delta^l$, meets a mark of  $\{\widehat{E}_t^{i,+,
  \widehat{m}}(y+\Delta^l)+\widehat{E}_t^{i,-, \widehat{m}}(y+\Delta^l)\}$.  
Therefore, this discrepancy, corresponding to the stopping time 
$\tau_{l+1}$, is the first event of the \emph{$l$-th-discrepancy process} 
$\{D^l_t\}_{t \in [\tau_l, \infty)}$, given by 
\be
\begin{split}
D^l_t=&\sum_{i=1}^k \Bigl[\widehat{E}_t^{i,+, m}(y_t^l)+ \widehat{E}_t^{i,-, m}(y_t^l)\\
&\qquad{}+\widehat{E}_t^{i,+,\widehat{m}}(y_t^l+\Delta^l)+\widehat{E}_t^{i,-,
  \widehat{m}}(y_t^l+\Delta^l)\Bigr]
\end{split}
\ee
where $y_t^l$ is  the density profile 
path  defined by a realization of the (level $l$) construction 
done at this stage.
\medskip
 
The construction done at the $l$-th stage makes sense, 
and has the correct rates, for
$t \geq \tau_l$.  Thus, together with the assumed graphical construction 
for $t\in[0,\tau_l)$, it yields a well defined coupling for
$g_t^{x^0,N}$ and $\widehat{g}_t^{x^0,N}$ at all
times.  Such a (level-$l$) coupling, however, loses precision after the
next discrepancy is encountered.  To improve it, we ignore
it for $t\geq \tau_{l+1}$ and replace it by the level-$(l+1)$ construction
corresponding to the $l+1$ stage.  This stage begins with 
$g_{\tau_{l+1}}^{x^0,N}-\widehat{g}_{\tau_{l+1}}^{x^0,N}=\Delta_{l+1}$ 
with $|\Delta_{l+1}-\Delta_l|=1/N$.

This recursive construction is continued, for each trajectory, until the time
$t=NT$ is achieved.  The procedure involves, almost surely,
 a finite number of stages.  The process 
 \be
 \label{er:dd}
 \overline D_t \;=\; D^l_t  \quad \mbox{if } t\in(\tau_l,\tau_{l+1}]
 \ee
 $l=0,1,\dots$ ($\tau_0=0$), counts the number of discrepancies.  
 It satisfies the relation $\{\overline{D}_t\geq l\}=\{\tau_l \leq t\}$.

\subsection{Discrepancy rates}
\label{sr.cd}

The proof of Theorem \ref{theo} requires the control of the
distance between $m_t^{x^0,N}$ and $\widehat{m}_t^{x^0,N}$.  As each
discrepancy brings an additional separation of $1/N$, 
\be
\label{er:mm}
\bigl| m_t^{x^0,N} - \widehat{m}_t^{x^0,N}\bigr| \;\le\; \frac{\overline D_{Nt}}{N}
\ee
To estimate the right-hand side we 
first determine upper bounds on the time-dependent rate of the process
$\{\overline D_t\}$.

\begin{lemma}
\label{lr.rdisc}
Consider $N\in\mathbb{N}$, $T\ge 0$ and $\delta>0$.
For each $l\in\mathbb{N}$, let $R^l_t$ be the instantaneous rate of the
level-$l$ discrepancy process $D_t^l$, $t\in[\tau_l,\tau_{l+1}]$ defined above
and let $R^l=\sup\bigl\{R^l_t: t\in[\tau_l,\tau_{l+1}]\cap [0,NT]\bigr\}$.  Then 
there exists a constant $A>0$ such that the events
\be
\label{er:rnt}
\mathcal{R}^{NT}_\delta \;=\; \Bigl\{ R_l \le N^{\delta -1/2} + \frac {A\,l}{N} \quad \forall \,l
\mbox{ s.t. } \tau_l\le NT\Bigr\}
\ee
satisfy
\be
\label{er:lii}
P\bigl(\underline \lim_N \mathcal{R}^{NT}_\delta \bigr) \;=\; 1\;.
\ee
\end{lemma}

\proof
Let $\Delta_t$ be the distance between the coupled geometrical realizations 
$g^{x^0,N}_t$ and $\widehat g^{x^0,N}_t$: 
\be
\Delta_t=\sum_{l\geq 0} \Delta^l {\bf 1}_{t \in [\tau_l,\tau_{l+1})}.
\ee
The discrepancy process can be written as
\be
\begin{split}
\overline{D}_t=&\sum_{i=1}^k\Bigl[
\widehat{E}_t^{i,+, m}(g^{x^0,N}_t)+ \widehat{E}_t^{i,-,
  m}(g^{x^0,N}_t)\\ 
&\qquad {}+\widehat{E}_t^{i,+,\widehat{m}}(g^{x^0,N}_t+\Delta_t)+\widehat{E}_t^{i,-,
  \widehat{m}}(g^{x^0,N}_t+\Delta_t)\Bigr]\;.
\end{split}
\ee
The rate of this process is zero at $t=0$, but it increases
as the processes $\bigl\{g^{x^0,N_t}\bigr\}$, $\bigl\{\widehat g^{x^0,N}_t\bigr\}$
and $\bigl\{x_t^{x^0}\bigr\}$ move away
from each other during the stochastic evolution.
 
We see from (\ref{+}) and (\ref{-}) that to bound this rate 
we must compare values of $x\,\lambda(y)$ and $x\,\mu(y)$ for
different densities $x$ and $y$. Due to the
Lipschitz hypothesis (\ref{lipchitz}), these difference increase at most
linearly, and there exists a constant $A$
such that the rate of $\overline{D}_t$ is bounded above by 
\be
\label{x}
 A\,\bigl|\widehat{g}_{t}^{x^0,N}-x_{t/N}^{x^0}\bigr|+
 A\,\bigl|g_{t}^{x^0,N}-\widehat{g}_{t}^{x^0,N}\bigr|\;.
\ee

For a given realization of the graphical construction,
the second term in (\ref{x}) is bounded above by 
$\frac{\overline{D}_t}{N}$, as remarked in (\ref{er:mm}).  Therefore
\be
\label{er:addl}
\bigl|g_{t}^{x^0,N}-\widehat{g}_{t}^{x^0,N}\bigr| \;\le\; \frac lN
\qquad \mbox{if } t\in[\tau_l,\tau_{l+1}]\;.
\ee

For  the first term in (\ref{x}) we apply first the probabilistic bound 
\begin{eqnarray}
\label{local}
\lefteqn{
P\Bigl(\bigl|\widehat{g}_{t}^{x^0,N}-x_{t/N}^{x^0}\bigr|> \frac{1}{2A} N^{\delta-1/2}\Bigr)}
\nonumber\\
&&=\;P\Bigl(\bigl|\widehat{m}_{t/N}^{x^0,N}-x_{t/N}^{x^0}\bigr|> \frac{1}{2A} N^{\delta-1/2}\Bigr)\\
 &&< \;\exp (-c N^{\delta})\;.\nonumber
\end{eqnarray}
valid for \emph{each} $t>0$.  The last inequality follows from  (\ref{bound2}).
We need, however, a bound valid for \emph{all} $t\in[0, NT]$.  To this end, we
apply (\ref{local}) to a sufficiently thick collection of times.  We pick a positive real $\gamma$
(soon to be chosen larger than 3) and denote $M$ the integer part of $N^\gamma$.
For each $0\leq j \leq M$ let
\be
C_j\;=\;\Bigl\{\bigl|\widehat{g}_{jNT/M}^{x^0,N}-x_{jT/(NM)}^{x^0}|\leq
N^{\delta-1/2}\Bigr\}
\ee
and
\be
\Theta \;=\; \inf \Bigl\{t: \bigl|\widehat{g}_{t}^{b(x^0),N}-x_{t/N}^{x^0}\bigr|>
N^{\delta-1/2}\Bigr\}\;.
\ee
Then, 
\be
\begin{split}
\label{global}
P(\Theta \leq  NT) & \leq\;P\bigl(\Theta \leq  NT, \cap_{l=0}^{M} C_l\bigr) + 
\sum_{l=0}^M \bigl(1-P(C_l)\bigr)  \\
&\leq\;M\Bigr[1-\Bigl(1-\frac{dNT}{M}\Bigr)\exp(-dNT/M)\Bigr]+ 
M\exp (-c N^{\delta}) \\
& \leq c\, N^{2-\gamma}
\end{split}
\ee
where $d$ and $c$ are positive constants. In the second line we used
(\ref{local}) to bound the last term. For the other term, we just observed
that the conditions $\Theta \leq  NT$ and $\cap_{l=0}^{M}
C_l$ together imply that the process must have at least two transitions during
the time interval of length $NT/M$ containing $\Theta$. 
The constant $d$ bounds the rate of flips of the process 
$\{\widehat{g}_{t}^{b(x^0),N}\}_{t\geq 0}$
(we can take
$d= \sum_{i=1}^k \|\lambda_i\|_\infty+  \sum_{i=1}^k \|\mu_i\|_ \infty$).  
The choice $\gamma >3$ yields a summable 
bound in (\ref{global}), which implies
\be
\label{asresult}
P\Bigl(\overline{\lim}_{N} \bigl\{\Theta \leq  NT\bigr\}\Bigr)=0
\ee
This result together with the bound (\ref{er:addl}) proves (\ref{er:lii}). \qed

\subsection{Conclusion of the proof}

Due to Lemma \ref{lr.1} and relation (\ref{er:mm}), the following lemma 
concludes the proof of Therem \ref{theo}.

\begin{lemma} 
For any $\varepsilon >0$ and $0 \leq t\leq T$,
\be
\label{nt}
P\Bigl(\overline{\lim}_{N} \bigl\{\overline{D}_{NT} \geq 
N^{\varepsilon +1/2 }\bigr\}\Bigr)\;=\;0\;.
\ee
\end{lemma}

\proof
Let us denote $\widetilde{N}_t= \overline{D}_{Nt}$.  This process 
---which has rates $N$ times
higher than those of $\{\overline{D}_t\}_{t\geq0}$--- counts
discrepancies in the time scale of $\{m_t^{x^0,N}\}_{t\geq 0}$. 
Let $\mathbf T_N$ be the time needed for the latter to collect 
$N^{\varepsilon + 1/2}$ discrepancies:  
\be
\mathbf T_N \;=\; \min\bigl\{ t: \widetilde N_t \ge N^{\varepsilon + 1/2}\bigr\}\;.
\ee
It can be written in the form
\be
\label{er:prt}
T_N= \sum_{i=1}^{N^{\varepsilon + 1/2}_+}T_i
\ee
where $T_1, T_2, \ldots $ are the independent successive times spent in
between jumps and $N^{\varepsilon + 1/2}_+$ is the smallest integer 
following $N^{\varepsilon + 1/2}$.

We now choose some $\delta$ with $0<\delta<\varepsilon$.  By Lemma \ref{lr.rdisc}, the events 
\be
\mathcal D_r=\bigl\{\mbox{condition (\ref{er:rnt}) is valid for } N\ge r\bigr\}
\ee
satisfy  
\be
P\Bigl( \bigcup_{r\in\mathbb N} \mathcal D_r \Bigr) \;=\; 1\;.
\ee
In the sequel we shall show that 
\be
\label{er:sss}
\sum_N P\bigl(\mathbf T_N<T \,;\, \mathcal D_r\bigr) \;<\; \infty
\ee
for each natural number $r$.  
This concludes the proof because it implies that
\be
P\Bigl(\overline{\lim}_N\bigl\{\mathbf T_N<T\bigr\}\Bigr) \;\le\;
\sum_r P\Bigl(\overline{\lim}_N\bigl\{\mathbf T_N<T\bigr\}\,;\, \mathcal D_r\Bigr)
\;=\;0\;.
\ee

To prove (\ref{er:sss}) we partially resum the decomposition (\ref{er:prt})
in blocks of size 
\be
\label{er:qqq}
Q\;=\;\frac{N^{\varepsilon + 1/2}_+}{N^{\delta + 1/2}_+}
\;\sim\; N^{\varepsilon-\delta} \;\tende{N\to\infty} \;\infty
\ee
We consider intervals
$I_l=[(l+1)N^{\delta+1/2}_+ +1,lN^{\delta+1/2}_+]$ and write 
\be
{\bf T}= \sum_{l=1}^{Q} G_l \quad,\quad G_l=\sum_{j\in I_l} T_j\;.
\ee
Within $\mathcal D_r$ the process $\{\tilde{N}_t\}_{t\geq0}$  jumps 
from $i$ to $i+1$ with rates bounded above by $N^{\delta+1/2} +
A\,i$. Thus, for each $i \in I_l$ the rate of $T_i$ is bounded above by
$N^{\delta+1/2} + A\,l$, which is smaller than 
$(l+1)N^{\delta+1/2}$ if $N$ is large enough. 
This shows that, for such $N$'s, the output of each variable $G_l$
is no smaller than that of a sum of  
$N^{\delta+1/2}_+$ i.i.d. exponential random variables 
with rate $(l+1)N^{\delta+1/2}_+$.  Hence, 
\be
P\bigl({\bf T} < T\,;\, \mathcal D_r\bigr ) \;\leq\; 
P\biggl( \sum_{l=1}^{Q} \frac{G_l(N^{\delta+1/2}_+)}{(l+1)N^{\delta+1/2}_+} <T\biggr)
\ee
where $\{G(N^{\delta+1/2})\}_{l \geq 1}$
denotes an i.i.d. sequence of Gamma random variables with
with parameters $n=N^{\delta+1/2}_+$ and $\lambda=1$. 
The large-deviation properties of these distributions imply that
each variable $G_{l,N}=G_l(N^{\delta+1/2})/N^{\delta+1/2}$ satisfies 
\be
\label{er:lade}
P( G_{l,N} <1/2) \;\leq\;
\exp(-cN^{\delta+1/2}) 
\ee
for some positive constant $c$ and all $1\leq l \leq Q$ and $N$ large enough.
\smallskip

Denoting $A_{l,N}=\{G_{l,N} \geq 1/2)\}\cap\mathcal D_r$
and $B_{Q,N}=\cap_{l=1}^{Q} A_{l,N}\cap\mathcal D_r$ we have
\be
\label{t}
P\bigl({\bf T}_N < T\,;\,\mathcal D_r\bigr) \;\leq\; \bigl(1-P(B_{Q,N})\bigr) 
+ P\bigl({\bf T}_N<T, B_{Q,N}\bigr)
\ee
On the event $B_{Q,N}$, ${\bf T}$ is bounded below by 
\be
\frac{1}{2}\sum_{l=1}^Q \frac{1}{l+1} \;\sim\; \log Q
\;\tende{N\to\infty}\;\infty\;.
\ee
Therefore the second term in the right-hand side of (\ref{t})
 is zero for $N$ large enough.  Bounding the first term by the large-deviation
 estimate (\ref{er:lade}) we conclude that 
 \be
P\bigl({\bf T}_N < T\,;\,\mathcal D_r\bigr) \;\leq\; 
Q \, \exp(-cN^{\delta+1/2})
\ee
for $N$ large enough.  This proves (\ref{er:sss}). \qed

\section*{Acknowledgments}
 The authors wish to thank the NUMEC and the University of S\~ao Paulo (R.F.)   
 and the University of Rouen (E.J.N.) for hospitality during 
 the completion of this work and to the USP-COFECUB agreement. 
E.J.N.~also thanks Antonio Prudente Cancer Research (FAPESP-CEPID).
L.R.F.~is partially supported by the CNPq grant 307978/2004-4.
E.J.N.~and L.R.F.~are partially supported by the by FAPESP grant 2004/07276-2. 
All three authors are partially suported by CNPq grant
484351/2006-0. R.F.~benefited from the CNRS-FAPESP agreement.
 

\end{document}